
\baselineskip=14pt
\parskip=10pt

\font\eightrm=cmr8 
\font\eighttt=cmtt8
\magnification=\magstephalf
\def\F{{\cal F}}

\def\1{{\overline{1}}}
\def\2{{\overline{2}}}
\parindent=0pt
\overfullrule=0in

\def\frac#1#2{{#1 \over #2}}
\bf
\centerline
{
Automatic Proofs  of Asymptotic ABNORMALITY (and much more!) of  Natural Statistics 
}
\centerline
{
Defined on Catalan-Counted Combinatorial Families
}
\rm
\bigskip
\centerline{ {\it
Shalosh B. EKHAD and
Doron 
ZEILBERGER}\footnote{$^1$}
{\eightrm  \raggedright
Department of Mathematics, Rutgers University (New Brunswick),
Hill Center-Busch Campus, 110 Frelinghuysen Rd., Piscataway,
NJ 08854-8019, USA.
{\eighttt zeilberg  at math dot rutgers dot edu} ,
\hfill \break
{\eighttt http://www.math.rutgers.edu/\~{}zeilberg/} .
March 21, 2014.
Accompanied by the Maple packages \hfill\break
{\eighttt http://www.math.rutgers.edu/\~{}zeilberg/tokhniot/AlgFunEq }, and  \hfill\break
{\eighttt http://www.math.rutgers.edu/\~{}zeilberg/tokhniot/Cheyne } \quad . \hfill\break
Supported in part by the NSF. Exclusively published in the Personal Journal of Shalosh B. Ekhad and
Doron Zeilberger, and arxiv.org .
}
}

{\bf Preliminary Sermon: Humans will be Humans; The Medium is the Message}

The famous  Catalan numbers (see [Sl1]), count zillions of
combinatorial families (see [St])
and many humans have fun trying to find `nice' bijections between family A and family B.
While this may be {\it fun} for a while, sooner or later this game gets old, especially since
the {\it real reason} Catalan numbers are so ubiquitous is their {\it simplicity},
and that humans can {\it only} grasp simple things.

Indeed, (see [Z]), the reason for the {\it ubiquity} of the 
sequence of Catalan numbers, $\{ c_n \}$, is that their
generating function
$$
C(z):= \, \sum_{n=0}^{\infty} c_n z^n \quad ,
$$
satisfies the {\bf simplest} possible (genuinely!) {\it algebraic} equation, namely
$$
C(z)=1+zC(z)^2 \quad ,
$$
that is equivalent to the quadratic recurrence satisfied by the Catalan numbers themselves,
namely:
$$
c_n \, = \, \sum_{k=1}^{n} c_{k-1} c_{n-k} \quad , \quad c_0=1 \quad .
$$

Often, the members of the combinatorial family in question posses natural {\it statistics},
for example for Dyck paths, the number of `inversions' ($D$ (not necessarily immediately) ahead of $U$),
or for $132$-avoiding permutations, the number of occurrences of some given pattern, then
it may happen that two different statistics `amazingly' have the same average!
{\bf Wow!},  Let's find a bijection! See, e.g., the humanly-generated  article [B] (by human Mikl\'os B\'ona),
that appeared in the very prestigious (and very selective!) {\it Electronic journal of Combinatorics}, 
that does its best [alas, not always successfully] to {\it only} accept {\it the best papers}.
[It often errs on both sides, rejecting truly seminal papers, and accepting quite a few trivial ones.]

Humans can, with some effort, find closed-form expressions for the {\it average} (aka {\it expectation}, aka {\it first moment}),
of a given combinatorial {\it statistic} (aka {\it random variable}),
and if they try {\bf really} hard, {\it may} be able to find the {\it variance} (aka {\it second moment} [about the mean]), 
but beyond that they should enlist their much superior silicon brethrern, and 
develop algorithms  for discovering (and proving!) closed-form expressions for as many as possible  moments. 
In addition to its intrinsic interest, this activity
would also indicate whether the combinatorial
statistic in question seems to be {\it asymptotically normal} (if the standardized moments, starting
with the third, converge, as $n$ goes to $\infty$,  to $0,3,0,15,0,105, \dots$, the famous moments of the normal distribution), or
whether it is (rigorously-) {\bf provably} {\it not} normal 
(if the expression for the {\it skewness} (alias {\it standardized 3rd moment}) does not
tend to zero, we are done!) .

In the present article, a collaboration between a human (DZ) and computer (SBE) we do just that! The
human wrote a Maple package available, free of charge from:

{\tt http://www.math.rutgers.edu/\~{}zeilberg/tokhniot/AlgFunEq } \quad ,

that {\it once written}, can handle {\bf zillions} of possible statistics defined on Catalan-families,
and {\it surprise-surprise}, find {\bf zillions} of B\'ona-style `surprises', and of course,
prove them all fully rigorously.
More importantly, it can prove  {\bf asymptotic abnormality}, by deriving (and proving!)
closed-form expressions for the  skewness, as an expression in $n$, and having done that,
(automatically!) take the limit as $n \rightarrow \infty$, and realize that
that limit is {\bf not} zero. Since for any asymptotically normal sequence of random variables,
that limit should be zero, this constitutes a fully rigorous proof of asymptotic abnormality.
But why stop with the skewness? Our program also finds closed-form expressions for the {\it kurtosis},
aka {\it standardized fourth moment} (and proves that its limit, as $n \rightarrow \infty$, is {\it not} $3$), as well
as expressions for  higher moments.

It is true that, with great effort, very smart humans, like
Svante Janson ([J]), can do it by entirely human means (and can even handle {\it all} moments,
at least recursively), but they can {\bf only} do the leading asymptotics!
Not even Svante Janson can find, {\it just by hand}, e.g., an {\bf exact} closed-form expression,
in $n$, for the sixth moment of the random variable `number of occurrences of the pattern $213$'
in the set of $132$-avoiding permutations of length $n$.

But do we really care about the 6th moment of some stupid statistic defined on some stupid 
family of sets?
Of course not! {\bf The Medium is the Message!} This article is but a {\it case-study}
in {\bf human-computer collaboration}. The human teaching the computer how to
solve every conceivable problem in a wide class of  combinatorial  problems,
by the human {\it designing} algorithms, then {\it implementing} them
(in our case in Maple),
and then letting the computer {\bf execute} them. Once we get better and better at this kind of collaboration,
we would be ready for the big time! Stand by (in 100 years or less) for a computer-generated
proof of RH and $NP \neq P$.

{\bf Maple packages and Sample Input and Output Files}

As usual in the ongoing collaboration between the authors of the present article,
the most important part is {\bf not} the {\it article}, but the {\bf Maple packages} that come with it,
that can be dowloaded, {\it free of charge}, from the front of this article

{\tt http://www.math.rutgers.edu/\~{}zeilberg/mamarim/mamarimhtml/abnormal.html} \quad ,

and the numerous input and output files also available there. In particular, the file

{\tt http://www.math.rutgers.edu/\~{}zeilberg/tokhniot/oAlgFunEq1BonaRedux} \quad ,

reproduces, in {\bf 0.117 seconds}, all the results  (all rigorously proved!)
of the first section of [B], that handled averages of the number
of occurrences of patterns of length $\leq 3$ (see below for details).

Before describing our new algorithms, we challenge
our readers (both humans and machines), with a {\bf neat conjecture},
based on ample evidence outputted by our Maple package {\tt AlgFunEq}.

{\bf Lots and Lots of B\'ona-Style Surprises and a Conjecture}

Let $AV_n(132)$ be the set of $132$-avoiding permutations of $\{1 , \dots , n\}$,
and for any permutation $\pi$ and pattern $p$, let
$a_p(\pi)$ be the number of occurrences of the pattern $p$ in the
permutation $\pi$ (in other words, if $\pi$ has length $n$ and $p$
has length $k$, the number of $k$-tuples 
$$
1 \leq i_1 < i_2 < \dots < i_k \leq n \quad,
$$
such that $\pi_{i_1} \pi_{i_2} \dots \pi_{i_k}$ is `order-isomorphic' to $p$).
For each pattern $p$, define the sequence
$$
A_p(n):=\sum_{\pi \in AV_n(132)} a_p(\pi) \quad .
$$

In [B], B\'ona observed that trivially $A_{231}(n)=A_{312}(n)$
(since $231$ and $312$ are inverses of each other, and the
class of $132$ permutations is closed under taking inverse),
but that {\it surprisingly} (at least to him),
both are also equal to $A_{213}(n)$.
Hence we have the following facts.

$\bullet$ For $k=1$ there is $1$ B\'ona class (of course!)

$\bullet$ For $k=2$ there are $2$ B\'ona classes (of course!)

$\bullet$ For $k=3$ there are $3$ B\'ona classes 

Can you spot a {\it pattern}?(pun intended!). Hint: the term for $k=4$ is {\bf not} $4$.

The first-named author

(See {\tt http://www.math.rutgers.edu/\~{}zeilberg/tokhniot/oAlgFunEq1} \quad , \hfill\break
{\tt http://www.math.rutgers.edu/\~{}zeilberg/tokhniot/oAlgFunEq1a9} \quad and \hfill\break
{\tt http://www.math.rutgers.edu/\~{}zeilberg/tokhniot/oAlgFunEq1a10})

{\it rigorously} proved that the numbers of distinct {\it B\'ona classes} (i.e. distinct sequences $\{A_p(n)\}$
as $p$ ranges over all $132$-avoiding patterns of length $k$) for $1 \leq k \leq 10$,
starting with $k=1$, are as follows:
$$
1, 2, 3, 5, 7, 11, 15, 22, 30, 42,  \dots \quad .
$$
Going to the {\it indispensable} OEIS, 
we immediately
realized that these are the first ten values of [Sl2],
and this naturally leads to the following  intriguing conjecture.

{\bf Conjecture} ($100$ donation to the OEIS in honor of the prover or disprover):
For every $k \geq 1$, the number of {\it distinct} sequences
$A_p(n)$, as  $p$ ranges over all the $(2k)!/(k!(k+1)!)$ $132$-avoiding patterns
of length $k$, is exactly $p(k)$, the number of integer partitions of $k$.

[Of course the letter $p$ in `$p(k)$' has no relation whatsoever to the letter $p$ in ``$A_p(n)$'',
except that both words, `partition' and `pattern', happen to start with it.]

Ideally one would like to have not just explicit expressions for the {\it average} (expectation)
of the random variable `number of occurrences of the pattern $p$',
i.e. $A_p(n)/c_n$, {\bf but} explicit expressions (or failing this,
{\it efficient} algorithms for generating many terms) for computing as many as possible {\it higher moments}.

{\bf Higher Moments}

Every infinite sequence of sets,  let's call it $\{ C_n \}_{n=0}^{\infty}$, counted by the
Catalan numbers, i.e. such that $|C_n|=c_n$ is (most probably) so because of a {\it natural} structure-bijection
[sometimes obvious (e.g. binary trees, Dyck paths), sometimes less so (e.g. $123$-avoiding permutations)]
$$
C_n \leftrightarrow \bigcup_{k=1}^{n} C_{k-1} \times C_{n-k} \quad, 
$$
leading immediately to the famous recurrence
$$
c_n = \sum_{k=1}^{n} c_{k-1} c_{n-k} \quad , \quad c_0=1 \quad .
$$

Many times, the members of our Catalan family have
{\it personalities}, and posses {\it numerical attributes} (usually positive integers, but
not necessarily), interchangeably called {\it statistics} and {\it random variables}.
Let's call such a statistic $s \rightarrow i(s)$. Then a natural question is

$\bullet$ What is the {\it average}, getting a brand-new numerical sequence

$$
a_n := \frac{1}{c_n} \sum_{s \in C_n} i(s) \quad .
$$

But why stop here? For each power $r$, we may be interested in the $r$-th moment, getting yet another numerical
sequence, one for each $r$,

$$
m^{(r)}_n := \frac{1}{c_n} \sum_{s \in C_n} (i(s))^r \quad .
$$

For statisticians (and probabilists),  more insightful sequences are {\it moments about the mean},

$$
M^{(r)}_n := \frac{1}{c_n} \sum_{s \in C_n} (i(s)-a_n) ^r \quad ,
$$

that of course, are easily computed from  $m^{(r)}_n$, using the binomial theorem.
For $r>2$, in fact, the most interesting quantities are the {\it standardized moments},
aka {\it alpha coefficients}, 
$$
\alpha^{(r)}_n= \frac{M^{(r)}_n}{ (M^{(2)}_n)^{r/2} } \quad,
$$
that almost always converge, as $n \rightarrow \infty$, to a sequence of real numbers, let's call them $\beta^{(r)}$.
When that happens, there is a {\it limiting distribution}, that in many cases (but not for Catalan families!)
is the good-old {\it Gaussian} (aka {\it normal}) distribution, and that happens when
$\beta^{(r)}=0$ for $r$ odd and $\beta^{(r)}=1 \cdot 3 \cdots (r-1)=r!/(2^{r/2}(r/2)!)$, for $r$ even.

{\bf Weighted-Counting of Catalan Objects According to a Statistic}

Knowing all the moments is equivalent to knowing {\it explicitly} the {\it generating function}
(aka {\it weight-enumerator}) according to the statistic $i(s)$, using the indeterminate $t$,
getting a family of {\it polynomials}

$$
P_n(t) := \sum_{s \in C_n} t^{i(s)} \quad .
$$

Once we know $P_n(t)$ we can expand it in terms of $t-1$
$$
P_n(t)=\sum_{r=0}^{\infty} \frac{f^{(r)}_n}{r!} (t-1)^r \quad,
$$
immediately getting 
$$
f^{(r)}_n := \sum_{s \in C_n} r! {{i(s)} \choose{r}}  \quad ,
$$
from which the {\it factorial moments} can be gotten upon dividing by $P_n(1)=c_n$.

From the factorial moments, the usual moments can be easily computed  (using Stirling numbers of the second kind).

Unlike the (numerical) enumerating sequence, $c_n$, that has a lovely {\it closed-form}, namely the famous
$$
c_n = \frac{(2n)!}{n!(n+1)!} \quad ,
$$
it is (usually) too much to hope for a closed-form expression for $P_n(t)$, and in fact, their generating function
is usually not even algebraic, i.e. the `grand-generating function', w.r.t to $z$, say
$$
\F(z,t):=\sum_{n=0}^{\infty} P_n(t) z^n \quad ,
$$
usually does {\it not} satisfy an analogous algebraic equation to $C(z)=1+zC(z)^2$.

But it so happens (in many cases!), that the generating functions for the average (times $c_n$) and
the $r$-th factorial (and hence actual) moments (again times $c_n$), for each {\it specific} (i.e. numeric) $r$ are {\it algebraic}! In fact
the polynomial equations satisfied by those generating functions 
often happen to be of degree $2$, just like the one for the Catalan numbers, but, of course, with much more complicated
coefficients. How can me find them?

{\bf Functional Recurrences to the Rescue}

For purely {\it pedagogical} reasons, let's first consider a very simple example. 
As in [B] and [J], our {\it population} is the
set of $132$-avoiding permutations, i.e. the set of permutations, $\pi$, of $\{1, \dots , n\}$, such
that you {\bf never} have $1 \leq i_1 < i_2 < i_3 \leq n$ with  $\pi_{i_1} < \pi_{i_3} < \pi_{i_2}$.
Let's first convince ourselves that this is indeed a Catalan family.

Take a typical such permutation, $\pi$, and look for the {\it location} of the largest entry, $n$.
Suppose $n$ stands at the $k$-th place, i.e. $\pi_k=n$. Then it is easy to see that
{\bf all} the  entries standing to the {\it left} of $n$, i.e. $\{ \pi_1, \dots , \pi_{k-1} \}$ are
all {\bf larger} than all the entries standing to the right of $n$, namely $\{ \pi_{k+1}, \dots, \pi_n \}$,
or else a forbidden $132$ pattern would emerge with the $n$ playing the role of the `$3$' in $132$.

Hence every such permutation can be written as
$$
\pi= \pi_1 n \pi_2 \quad ,
$$
where $\pi_1$ is a permutation of the set $\{n-k+1, \dots, n-1\}$ and $\pi_2$ is a permutation of
$\{1, 2, \dots, n-k \}$. Of course, both $\pi_1$ and $\pi_2$ are $132$-avoiding on their own right,
and the map is a bijection. 
Hence the number of $132$-avoiding permutations of $\{1, \dots, n\}$ with $\pi_k=n$ equals
$c_{k-1} c_{n-k}$, and summing over $1 \leq k \leq n$ yields the Catalan recurrence,
$c_n=\sum_{k=1}^{n} c_{k-1} c_{n-k}$, for the cardinality of the set of $132$-avoiding permutations.

But now let's consider the simple statistic `number of $21$ patterns'.

Let $a_{21}(\pi)$ be the number of $21$ patterns of $\pi$. Using the above decomposition $\pi= \pi_1 n \pi_2$, we clearly have
$$
a_{21}( \pi_1 n \pi_2)= a_{21}(\pi_1) + a_{21}(\pi_2) + k(n-k) \quad,
$$
since a $21$ pattern may either be entirely contained in $\pi_1$, entirely contained in $\pi_2$, or the `$2$' may
belong to $\pi_1$ ($k-1$ possibilities) and the `$1$' may belong to $\pi_2$ ($n-k$ possibilities),
so altogether $(k-1)(n-k)$ possibilities, and of course, the $2$ may be  the `$n$', and  that gives $n-k$ extra
scenarios, so altogether we have $(k-1)(n-k)+(n-k)=k(n-k)$ additional occurrences of the patter $21$.

Let's define the {\it weight-enumerator},
$$
P_n(t) :=\sum_{\pi \in AV_{132}(n)} t^{a_{21}(\pi)} \quad .
$$
For $1 \leq k \leq n$, let $AV^{(k)}_{132}(n)$ be the subset of $AV_{132}(n)$ for which $\pi_k=n$, then of course
$$
P_n(t) :=\sum_{\pi \in AV_{132}(n)} t^{a_{21}(\pi)} 
\,\, =\,\, \sum_{k=1}^{n} \,\, \, \sum_{\pi \in AV^{(k)}_{132}(n)} t^{a_{21}(\pi)} 
$$
$$
=\sum_{k=1}^{n} \,\,\,
\sum_{ { {\pi_1 \in AV_{132}(k-1)} \atop {\pi_2 \in AV_{132}(n-k)} } } 
t^{a_{21}(\pi_1)+a_{21}(\pi_2) + k(n-k)} 
=\sum_{k=1}^{n}  t^{k(n-k)} \sum_{ { {\pi_1 \in AV_{132}(k-1)} \atop {\pi_2 \in AV_{132}(n-k)} } } 
t^{a_{21}(\pi_1)} t^{a_{21}(\pi_2)}
$$
$$
=\sum_{k=1}^{n}  t^{k(n-k)} \left ( \sum_{\pi_1 \in AV_{132}(k-1)} t^{a_{21}(\pi_1)} \right )
\left ( \sum_{\pi_2 \in AV_{132}(n-k)} t^{a_{21}(\pi_2)} \right )
$$
$$
=\sum_{k=1}^{n}  t^{k(n-k)} P_{k-1}(t) P_{n-k}(t) \quad ,
$$
and hooray!, we found the {\it  (non-linear) recurrence equation}
$$
P_n(t) \, = \, \sum_{k=1}^{n} t^{k(n-k)} P_{k-1}(t) P_{n-k}(t) \quad , \quad P_0(t)=1 \quad ,
\eqno(NLR)
$$
from which one can immediately get the first one hundred (or whatever) terms, but alas, no {\it closed form}.

Nevertheless, one can {\it easily} get explicit expressions for both the generating functions,
and the sequences themselves, for the average and higher moments,
$m^{(r)}_n$ (times $c_n$), for numeric $r$, up to any desired $r$. Of course as $r$ gets larger, the `explicit' expressions would get
more and more complicated, and there is (probably) no hope to get a {\it symbolic} expression
in $r$, but we do what we can.

One way to derive expressions for higher moments is {\it empirical}. After you crank out the first $100$ terms of the
sequence $P_n(t)$, {\it guess} (using, e.g. the built-in Maple package {\it gfun} developed by
Salvy and Zimmerman[SaZ]) explicit expressions for the numerical sequences $\{ P'_n(1) \}$, $\{ P''_n(1) \}$,
etc. But one can proceed purely `rigorously' as follows.

Suppose that we are only interested in the moments up to $r \leq R$, then write,
$$
P_n(1+z)=\sum_{r=0}^R \frac{1}{r!}f^{(r)}_n z^r + O(z^{R+1}) \quad .
$$
Now plug this in into the above non-linear recurrence $(NLR)$, with $t$ replaced by $1+z$
$$
\sum_{r=0}^R  \frac{1}{r!} f^{(r)}_n z^r + O(z^{r+1}) \quad .
= \sum_{k=1}^{n} (1+z)^{k(n-k)} \left ( \sum_{r=0}^R  \frac{1}{r!} f^{(r)}_{k-1} z^r + O(z^{R+1})) \right )
\left ( \sum_{s=0}^R  \frac{1}{s!}f^{(s)}_{n-k} z^s + O(z^{R+1})) \right ) \quad .
$$
Now use the binomial theorem to expand   $(1+z)^{k(n-k)}$ to order $R$:
$$
(1+z)^{k(n-k)}= 1+ k(n-k)z+ \dots + {{k(n-k)} \choose {R}}  z^R+ O(z^{R+1})  \quad,
$$
and compare the coefficients of $z^0,z, \dots, z^R$ on both sides, getting non-linear numerical recurrences.

Comparing the coefficient of $z^0$ we get the good-old recurrence for $c_n=f^{(0)}_n$.
Comparing the coefficients  of $z$ leads to the non-linear (numerical) recurrence
for the sequence $f^{(1)}_n$ that assumes that you already know $f^{(0)}_n$ (as indeed
you do, it being equal to the Catalan number $c_n$).
$$
f^{(1)}_n= \sum_{k=1}^{n}   k(n-k) f^{(0)}_{k-1} f^{(0)}_{n-k}+
\sum_{k=1}^{n}  f^{(1)}_{k-1} f^{(0)}_{n-k}+ \sum_{k=1}^{n}  f^{(0)}_{k-1} f^{(1)}_{n-k} \quad ,
$$
and so on and so forth.

From here we have a rigorous proof that {\it a priori}, the sequences $f^{(1)}_n$, $f^{(2)}_n$, etc.  all
have {\it algebraic} generating functions. One way is to teach the computer how to translate
the recurrences into a system of algebraic equations for the generating functions, and then
solve it, but a much more reasonable way is to use the non-linear recurrences (that now only involve
{\it numbers}) to crank out sufficiently many terms to guess the algebraic generating functions,
that can be justified {\it a posteriori} by plugging-in.

Analogously, for the pattern $12$ one has the non-linear recurrence
$$
P_n(t)= \sum_{k=1}^{n} t^{k-1} P_{k-1}(t) P_{n-k}(t) \quad , \quad P_0(t)=1 \quad ,
$$
(why?).

{\bf Patterns of length 3}

For the two patterns of length $2$, namely $12$ and $21$ we got away with simple (non-linear) recurrences,
but for patterns of length $3$, we have two new concepts. The first is (non-linear) 
{\it functional recurrence equation},
and the second one is {\it catalytic  variable}.

Let's try and find an analogous recurrence for the weight-enumerators
$$
P_n(t) :=\sum_{\pi \in AV_{132}(n)} t^{a_{231}(\pi)} \quad ,
$$
(Note that this new $P_n(t)$ is not the same as in the previous section, it is {\it local notation}).

So let's try to express 

$$
a_{231}(\pi)= a_{231}(\pi_1 n \pi_2) \quad ,
$$ 
in terms of $a_{231}(\pi_1)$  and $a_{231}(\pi_2)$, where 
$$
\pi_1 \in AV_{132}(k)[\{n-k+1, \dots,n-1\}] \quad ,
$$
and
$$
\pi_2 \in AV_{132}(n-k)[\{1, \dots,n-k\}] \quad .
$$
Here are the possible scenarios for an occurrence of the pattern $231$ in $\pi=\pi_1 n \pi_2$.

$\bullet$ it is completely immersed in $\pi_1$ : $a_{231}(\pi_1)$ ways.

$\bullet$ it is completely immersed in $\pi_2$ : $a_{231}(\pi_2)$ ways.

$\bullet$ the `$23$' part of the pattern $231$ belongs to $\pi_1$ and the `$1$' part belongs to $\pi_2$: $a_{12}(\pi_1)\cdot (n-k)$ ways.

$\bullet$ the `$2$' part of the pattern $231$ belongs to $\pi_1$, the `$3$' is $n$, and the `$1$' part belongs to $\pi_2$: $(k-1)(n-k)$ ways.

Hence
$$
a_{231}(\pi_1 n \pi_2) =  a_{231}(\pi_1) +  a_{231}(\pi_2) + a_{12}(\pi_1) \cdot (n-k) + (k-1)(n-k) \quad .
$$
Alas, we have an {\it uninvited} guest, $a_{12}(\pi)$, so we need to figure out how to express it in terms of
$a_{12}(\pi_1)$ and $a_{12}(\pi_2)$, but that's easy
$$
a_{12}(\pi_1 n \pi_2) =  a_{12}(\pi_1) +  a_{12}(\pi_2) + k-1 \quad .
$$
In order to figure out a recurrence for $P_n(t)$, we need to introduce a {\it catalytic variable}, $q$, that
takes care of the r.v. $a_{12}(\pi)$ and define:
$$
Q_n(t,q) :=\sum_{\pi \in AV_{132}(n)} t^{a_{231}(\pi)}  q^{a_{12}(\pi)} \quad.
$$
At the {\it end of the day}, once we have $Q_n(t,q)$, we would plug-in $q=1$ and get our desired $P_n(t)=Q_n(t,1)$, but
until then we would have to put-up with $q$. 

It would be convenient to define a {\it weight}
$$
Wt(\pi)(t,q):= t^{a_{231}(\pi)}  q^{a_{12}(\pi)} \quad.
$$
So we have
$$
Wt(\pi_1 n \pi_2)(t,q):= t^{a_{231}(\pi)}  q^{a_{12}(\pi)} 
= t^{  a_{231}(\pi_1) +  a_{231}(\pi_2) + a_{12}(\pi_1) \cdot (n-k) + (k-1)(n-k)}  q^{ a_{12}(\pi_1) +  a_{12}(\pi_2) + k-1 } 
$$
$$
=t^{(k-1)(n-k)} q^{k-1} \cdot
\, ( \, t^{  a_{231}(\pi_1) + a_{12}(\pi_1) \cdot (n-k)}  q^{ a_{12}(\pi_1)}  \, )
\cdot \,
( \, t^{ a_{231}(\pi_2)}  q^{ a_{12}(\pi_2)} \, )
$$
$$
=t^{(k-1)(n-k)} q^{k-1} \cdot \,
\, ( \, t^{  a_{231}(\pi_1)}  (qt^{n-k})^{a_{12}(\pi_1)} \, )
\, \cdot \,
\, ( \, t^{ a_{231}(\pi_2)}  q^{ a_{12}(\pi_2)}  \, )
$$
$$
=t^{(k-1)(n-k)} q^{k-1} \,\cdot \,  Wt(\pi_1)(t,t^{n-k}q) \, \cdot \,  Wt(\pi_2)(t,q) \quad ,
$$
leading to the {\it functional (non-linear) recurrence}
$$
Q_n(t,q)=\sum_{k=1}^{n} t^{(k-1)(n-k)} q^{k-1} \, Q_{k-1}(t, t^{n-k}q) \, Q_{n-k}(t,q) \quad, \quad Q_0(t,q)=1 \quad .
$$

Similarly, for the other four $132$-avoiding patterns, we have:

$\bullet$ $a_{123}(\pi)$ (with catalytic variable $q$ corresponding to $a_{12}(\pi)$):
$$
Q_n(t,q)=\sum_{k=1}^{n} q^{k-1} \, Q_{k-1}(t, tq) \, Q_{n-k}(t,q) \quad, \quad Q_0(t,q)=1 \quad .
$$

$\bullet$ $a_{321}(\pi)$ (with catalytic variable $q$ corresponding to $a_{21}(\pi)$):
$$
Q_n(t,q)=\sum_{k=1}^{n} q^{k(n-k)} \, Q_{k-1}(t, t^{n-k}q) \, Q_{n-k}(t,t^kq) \quad, \quad Q_0(t,q)=1 \quad .
$$

$\bullet$ $a_{213}(\pi)$ (with catalytic variable $q$ corresponding to $a_{21}(\pi)$):
$$
Q_n(t,q)=\sum_{k=1}^{n} q^{k(n-k)} \, Q_{k-1}(t, tq) \, Q_{n-k}(t,q) \quad,  \quad Q_0(t,q)=1 \quad .
$$

$\bullet$ $a_{312}(\pi)$ (with catalytic variable $q$ corresponding to $a_{12}(\pi)$):
$$
Q_n(t,q)=\sum_{k=1}^{n} q^{k-1} Q_{k-1}(t, q) Q_{n-k}(t,t^k q) \quad,  \quad Q_0(t,q)=1 \quad .
$$

All the above can be (and has been!) `taught' to the computer, and the computer can automatically derive
such functional equations. The above verbose derivation was only for the benefit of
explaining to humans the algorithms that would eventually be executed by computers.

Unfortunately, one can't get {\it closed-form} expressions for the $Q_n(t,q)$, and not even
a closed-form expression for their generating function, but 
the above functional recurrences are fairly efficient for
generating quite a few terms, and by plugging-in $q=1$ and taking successive derivatives
with respect to $t$ and then plugging-in $t=1$ one can generate quite large beginnings
of the moment-sequences, and have a guessing program (e.g. {\it gfun} [SaZ]) guess
either closed-forms or recurrences.

But for those who abhor guessing, one can get, completely automatically, algebraic
equations satisfied by the generating functions for the factorial moments, like
we did above for the $P_n(t)$ of $a_{21}(\pi)$. Now we need multi-variable Taylor
expansions, and need to put up with more multi-indexed sequences, but
so what? The computer does not mind!

{\bf How to represent Functional Recurrences in Maple?}

The beauty of mathematics, (and computers!) is that we can generalize and
consider a very general class of functional recurrences that include
all the above as very spacial cases.
\vfill\eject
Here goes:
$$
P_n(t_1, \dots,t_m)=
$$
$$
\sum_{k=1}^{n}
c(t_1, \dots, t_m,k,n)
P_{k-1} (a_1(t_1,\dots, t_m), \dots, a_m(t_1, \dots , t_m)) \cdot
$$
$$
P_{n-k}(b_1(t_1,\dots, t_m), \dots, b_m(t_1, \dots , t_m)) \quad , \quad P_0(t_1, \dots, t_m)=1 \quad .
$$

Here the $a_i$'s and $b_i$'s are {\it arbitrary} polynomials in their variables $t_1, \dots , t_m$, and
$c(t_1, \dots, t_m,k,n)$ is a polynomial in $t_1, \dots, t_m,k,m$ but where the powers of the $t_i$'s
(but not $n,k$) are allowed to be polynomials in $n,k$. For example
$c(t_1,t_2,k,n)=n^2t_1^{n+k^2}t_2^{n^3}+ k^3t_2^{k^3+n}$ is quite acceptable, but
$c(t_1,t_2,k,n)=k^k$ is {\bf not}.

In the Maple package {\tt AlgFunEq} such a Functional Equation is represented by the
following {\it data structure}
$$
[c, [a_1, \dots, a_m], [b_1, \dots, b_m]],[t_1, \dots, t_m] ] \quad .
$$

From this one can get, {\it automatically}, an efficient scheme for computing the mixed factorial moments
(and hence the pure ones). This is accomplished by procedure {\tt FAscheme}.

{\bf More General Functional Recurrences}

As general as the above form of the functional Catalan-type recurrences is, it is 
not general enough to consider many patterns of length larger than $3$. Let's take
for example $a_{4321}$.

We have
$$
a_{4321}(\pi_1 n \pi_2) =  a_{4321}(\pi_1) +  a_{4321}(\pi_2) + k a_{321}(\pi_2) + a_{21}(\pi_1) a_{21}(\pi_2) + (n-k) a_{321}(\pi_1) \quad.
$$
It is easy to see, because of the {\bf product} term  $a_{21}(\pi_1) a_{21}(\pi_2)$,
that an analogous derivation to the one for $a_{231}(\pi)$ carried above
does not lead to such a functional recurrence.

But what we {\it can} do is introduce ``generalized products''. First defining it  on pairs of monomials,
that the computer defines automatically, according to its needs,
and then extends it by bi-linearity to apply to any pair of polynomials. 
Calling this product $\F$, we have to handle functional equations of the 
more general form:
$$
P_n(t_1, \dots,t_m)=
$$
$$
\sum_{k=1}^{n}
c(t_1, \dots, t_m,k,n) \, \cdot \,
\F(P_{k-1} (a_1(t_1,\dots, t_m), \dots, a_m(t_1, \dots , t_m)) \, , \,
P_{n-k}(b_1(t_1,\dots, t_m), \dots, b_m(t_1, \dots , t_m))) \quad , 
$$
$$
P_0(t_1, \dots, t_m)=1 \quad .
$$
We confess that we were too lazy to implement these more general types in order to compute
higher moments, but for the special cases of {\it averages} (and that's all that B\'ona did for
patterns of length $3$) it is fully implemented for {\it any} pattern, see procedure {\tt MilonK}.

Note that even in this more general setting, we are guaranteed that the generating function
of each specific moment is algebraic, and hence it justifies the empirical guessing.

Using these, we were able to find closed-form expressions for the averages for {\it all} patterns
of length $\leq 10$, mentioned at the beginning of this article.

{\bf Rigorous Proofs of Asymptotic Abnormality}

Using the output of {\tt http://www.math.rutgers.edu/\~{}zeilberg/tokhniot/oAlgFunEq3} we see
(rigorously!) that the random variables $\pi \rightarrow a_p(\pi)$ are 
{\bf never} asymptotically normal, for patterns of length $\leq 3$,
in accordance with Janson's humanly-generated paper ([J]) (who proved it for patterns of all lengths).
That file also contains asymptotic expressions, to order $4$ of all standardized moments up to
the sixth.

{\bf What About The Number of Pattern-Occurrences in 123-Avoiding Permutations?}

The Catalan Structure of $123$-avoiding permutations is a bit more subtle.
For any $123$-avoiding permutation $\pi$, let $U(\pi)$ be the 
permutation obtained by finding the right-to-left maxima, circling them,
introducing an empty slot right before the first entry, sliding all the
non-left-to-right maxima one unit to the left, then increasing all the entries by $1$,
and finally sticking a $1$ at the remaining open slot. 

It is (almost) readily seen that any permutation $\pi$ of $\{1, \dots , n\}$ can be written, for some $k$, $1 \leq k \leq n$
$$
\pi= \pi_1 U(\pi_2) \quad ,
$$
where $\pi_1$ is a permutation of $\{k+1, \dots , n\}$ and $\pi_2$ is a permutation of $\{1, \dots , k-1\}$.

Cheyne Homberger[H] found generating functions, and explicit expressions, for  the {\it averages} of the random variables
``number of occurrences of pattern $p$'' defined on the set of $123$-avoiding permutations,
for all patterns $p$ of length $\leq 3$.

Let's be more general and try to find a functional recurrence equation for the weight-enumerator
$$
P_n(t):=\sum_{\pi  \in AV_{123}(n)} t^{a_{213}(\pi)} \quad .
$$
It turns out that we need two auxiliary statistics, and hence two catalytic variables:
$$
\sigma_1(\pi):=a_{213}(U(\pi))-a_{213}(\pi) \quad,
$$
$$
\sigma_2(\pi):=\sigma_1(U(\pi))-\sigma_1(\pi) \quad,
$$
and let's define the sequence of polynomials
$$
Q_n(t,s_1,s_2):=\sum_{\pi  \in AV_{123}(n)} t^{a_{213}(\pi)} s_1^{\sigma_1(\pi)}  s_2^{\sigma_2(\pi)} \quad ,
$$
then the reader is welcome to prove the (slightly) `weird' functional recurrence
$$
Q_n(t,s_1,s_2) \, = \, s_2 \, Q_{n-1}(t,s_1,s_2) \, +  \, \sum_{k=2}^{n} Q_{n-k}(t,s_1,s_2) \, Q_{k-1}(t,ts_1,s_1 s_2) \quad,
$$
subject to the initial condition $Q_0(t,s_1,s_2)=1$.
From this functional recurrence, one can compute quite a few terms. Then setting $s_1=1,s_2=1$, we get $P_n(t)=Q_n(t,1,1)$. See:

{\tt http://www.math.rutgers.edu/\~{}zeilberg/tokhniot/oCheyne2} \quad ,

where one can find the list of $P_n(t)$ for $0 \leq n \leq 35$.

Using this, we easily reproduced (a more compact version of) Homberger's explicit expression for what he called $a_n$ in [H], 
that is:
$$
\sum_{\pi \in AV_{123}(n)} a_{213}(\pi) \quad .
$$
The expression is: 
$$
-\frac{3}{8}\,{4}^{n}+\frac{1}{2}\, \left( n+2 \right)  \left( 2\,n-1 \right) c_{n-1} \quad,
$$
(recall that $c_n$ are the Catalan numbers).
Going beyond, we found an explicit expression for the {\bf second moment} (times $c_n$), i.e. for:
$$
\sum_{\pi \in AV_{123}(n)} (a_{213}(\pi))^2 \quad ,
$$
that turned out to be:
$$
-{\frac {9}{128}}\, \left( 3\,{n}^{2}+7\,n+6 \right) {4}^{n}+ \left( {\frac {19}{60}}\,{n}^{4}+{\frac {57}{20}}\,{n}^{3}+{\frac {67}{30}}\,{n}^{2}
+\frac{1}{10}\,n- 1 \right) c_{n-1} \quad ,
$$
but to our disappointment, the third moment turned out to be not nice at all, and all we could find was
a linear recurrence equation of order $4$ with coefficients that are polynomials of degree $4$, see the
output file

{\tt http://www.math.rutgers.edu/\~{}zeilberg/tokhniot/oCheyne3} \quad \quad .

All this was generated with the Maple package {\tt Cheyne} available from 

{\tt http://www.math.rutgers.edu/\~{}zeilberg/tokhniot/Cheyne} \quad \quad .

It should be be possible to expand this Maple package to handle larger patterns, but {\it enough is enough}.

{\bf Empirical Coda}

Recall that at the beginning of this article we conjectured that the number of `B\'ona classes' for $132$-avoiding permutations,
i.e. the number of {\it different} sequences that show up as averages of the random variables
`number of occurrences of the pattern $p$' for all $132$ patterns  of length $k$ is the number
of integer partitions of $k$, and we proved it {\it rigorously} for $k \leq 10$.

The output file

{\tt http://www.math.rutgers.edu/\~{}zeilberg/tokhniot/oCheyne1}

contains all the beginnings (through $n=9$) of the analogous sequences for $123$-avoiding permutations, 
and it is very possible that the sequence `total number of B\'ona classes' for 
patterns of length $k$ defined on the set of $123$-avoiding permutations, starting at
$k=1$ is:
$$
1,2,3,6,12,32, \dots  \quad \quad ,
$$
but we can't see a pattern. Can you?

{\bf References}

[B] Mikl\'os B\'ona, {\it Surprising symmetries in objects  counted by Catalan numbers},
The Electronic Journal of Combinatorics {\bf 19}(1) (2012), \#P62.
\hfill\break
{\tt http://www.combinatorics.org/ojs/index.php/eljc/article/download/v19i1p62/pdf}

[H]  Cheyne Homberger,  {\it Expected Patterns in Permutation Classes},
The Electronic Journal of Combinatorics {\bf 19}(3)(2012), \#P43
\hfill\break
{\tt http://www.combinatorics.org/ojs/index.php/eljc/article/download/v19i3p43/pdf}

[J] Svante Janson, {\it Patterns in random permutations avoiding the pattern 132}, arxiv.org, 22 Jan. 2014,
{\tt http://arxiv.org/abs/1401.5679}

[SaZ] Bruno Salvy and Paul Zimmerman, {\it GFUN: a Maple package for the manipulation of generating and holonomic functions in one variable},
ACM Transactions on Mathematical Software {\bf 20}(1994), 163-177 .

[Sl1] Neil Sloane, {\it Sequence \#108}, {\tt http://oeis.org/A000108} .

[Sl2] Neil Sloane, {\it Sequence \#41}, {\tt http://oeis.org/A000041} .

[St] Richard Stanley,  {\it Catalan Addendum}, {\tt http://www-math.mit.edu/\~{}rstan/ec/catadd.pdf} .

[Z] Doron Zeilberger, {\it Opinion \#49}, {\tt http://www.math.rutgers.edu/\~{}zeilberg/Opinion49.html}

\end